\documentclass[12pt,reqno]{amsart}
\usepackage{amssymb}

\newtheorem{theorem}{Theorem}[section]
\newtheorem{lemma}[theorem]{Lemma}
\newtheorem{corollary}[theorem]{Corollary}
\newtheorem{proposition}[theorem]{Proposition}

\theoremstyle{definition}

\newtheorem{definition}[theorem]{Definition}

\newtheorem{example}[theorem]{Example}

\newtheorem{remark}[theorem]{Remark}
\newtheorem{claim}[theorem]{Claim}

\def\ord{\operatorname{ord}}

\def\Z{\mathbb{Z}}
\def\Q{\mathbb{Q}}

\def\beq {\begin{equation}}
\def\endq {\end{equation}}
\def\bs{\boldsymbol}
\newcommand{\twolinesum}[2]{\sum_{\substack{{\scriptstyle #1}\\
{\scriptstyle #2}}}}

\newcommand{\Ocal}{\mathcal{O}}
\newcommand{\F}{\mathbb{F}}
\newcommand{\MOD}[1]{~\textup{mod}~#1}

\begin{document}
\title{Primitive divisors on twists of the Fermat cubic}
\thanks{The research of the second author is supported by a grant
from the NSERC of Canada.}
\author{Graham Everest, Patrick Ingram and Shaun Stevens}
\address{(GE+SS) School of Mathematics, University of East Anglia,
Norwich NR4 7TJ, UK}
\address{Department of Mathematics,
University of Toronto, Canada M5S 2E4} \email{g.everest@uea.ac.uk}
\email{pingram@math.utoronto.ca} \email{shaun.stevens@uea.ac.uk}

\begin{abstract}
We show that for an elliptic divisibility sequence on a twist of
the Fermat cubic, $u^3+v^3=m$, with $m$ cube-free, all the terms
beyond the first have a primitive divisor.
\end{abstract}

\maketitle


\section{Statement of Main Theorem}\label{intro}

Let $C$ denote a twist of the Fermat cubic,
\begin{equation}\label{homogform}
C: \quad U^3+V^3=mW^3
\end{equation}
with $m$ a non-zero rational number. If $K$ denotes any field of
characteristic zero, the set~$C(K)$ of projective $K$-rational
points satisfying (\ref{homogform}) forms an {\it elliptic curve}.
With respect to the usual chord and tangent addition the
set~$C(K)$ forms a group . The identity of the group is $(-1,1,0)$
and the inverse of the point $(U,V,W)$ is $(V,U,W)$. Let $R\in
C(\mathbb Q)$ denote a non-torsion rational point. Write
\begin{equation*}\label{defofUVW}nR=(U_n,V_n,W_n),\quad
U_n,V_n,W_n \in \mathbb Z
\end{equation*}
in lowest form with $\gcd(U_n,V_n,W_n)=1$. This paper
is devoted to proving the following theorem.

\begin{theorem}\label{theanswer} Let $C$ denote the elliptic curve
in~\eqref{homogform} with $m\in \mathbb Z$ assumed to be cube-free.
Let $W=(W_n)$ denote the sequence obtained as above from $R\in
C(\mathbb Q)$, a non-torsion rational point. For all $n>1$, the
term $W_n$ has a primitive divisor.
\end{theorem}

The term {\it primitive divisor} is defined in the following way.
\begin{definition}
Let $(A_n)$ denote a sequence with integer terms. We say an
integer $d>1$ is a {\it primitive divisor} of $A_n$ if
\begin{enumerate}\item $d\mid A_n$ and \item gcd$(d,A_m)=1$ for all non-zero
terms~$A_m$ with $m<n$.\end{enumerate}\end{definition}

The sequence $W=(W_n)$ is a {\it divisibility sequence}, which
means that, for all $m,n\in \mathbb N$,
\begin{equation}m\mid n \mbox{ implies } W_m\mid W_n.
\end{equation}
In line with recent developments~\cite{pe, emw, AMM, pi, pint, is,
sileds1, sileds2, sileds3} we define the sequence $W=(W_n)$ to be
an {\it Elliptic Divisibility Sequence}. Admittedly this stretches
the definition originally used by Morgan Ward~\cite{ward} but we
believe it is a reasonable name for a divisibility sequence that
arises from an elliptic curve.

In the sequel, it will often be convenient to work with the affine
curve
\begin{equation}\label{affinehomogform}
u^3+v^3=m
\end{equation}
and we will also refer to this curve as~$C$. Properties of points
such as being integral will generally refer to the affine curve.
Thus, an integral point on~$C$ is a pair of integers $(u,v)$
satisfying~(\ref{affinehomogform}). Theorem~\ref{theanswer} is the
best possible in the sense that if $R$ is an integral point, then
$W_1$ has no non-trivial divisors and so (by our definition) no
primitive divisors. Thus, the identity
$$(1+t)^3+(1-t)^3=(6t^2+2),$$
for any integer $t>1$, gives rise to an elliptic divisibility
sequence for which Theorem~\ref{theanswer} is best possible,
whenever $6t^2+2$ is cube-free. There are infinitely many~$t$ such
that $6t^2+2$ is cube-free by a result of Erd\H{o}s \cite{erdos}.

If one removes the condition that $m$ be cube-free, then it is
easy to construct counter-examples to Theorem~\ref{theanswer} by
clearing denominators. However, any elliptic curve in the form
(\ref{homogform}) can be transformed into one with $m$ a cube-free
integer by a simple scaling. Thus, Theorem~\ref{theanswer} does
give complete information, taking account of the transformation.

Note that Theorem~\ref{theanswer} has an immediate application to
the study of integral points on elliptic curves in the
form~(\ref{homogform}). In~\cite{taxi}, Silverman showed that,
when~$m$ is cube-free, there exists an absolute constant $\kappa$
such that~\eqref{homogform} has at most
$\kappa^{1+\textup{rank}(C/\Q)}$ integral points. In particular,
the number of integral solutions to~\eqref{affinehomogform} is
bounded uniformly if we restrict attention to curves $C$ with rank
not exceeding a given bound. Note that $C(\Q)$ is torsion-free
when its rank is positive. Thus if $\textup{rank}(C/\Q)\nobreak
=\nobreak 1$, then $C(\Q)$ consists solely of the multiples of a
single point $R$, say.  By Theorem~\ref{theanswer}, the
denominator of $nR$, for each $n\neq \pm1, 0$, has a primitive
divisor.  In particular, none of these denominators is 1.
\begin{corollary}\label{cor}
If $m\in \mathbb Z$ is cube-free, and $\textup{rank}(C/\Q)=1$,
then~\eqref{affinehomogform} has at most 2 integral solutions -
either of which generates $C(\Q)$.
\end{corollary}
Corollary~\ref{cor} is not overly surprising.  The method of proof
in \cite{taxi} gives very strong bounds on the number of integral
points on $C$ as the diophantine approximation involved in this
case is trivial. Corollary~\ref{cor} is noted because it is sharp
and completely qualitative. For quantitative results of similar
strength for other elliptic curves, see \cite{pint}.

In the next section, we will set Theorem~\ref{theanswer} in its
proper context, as well as outline the structure of the proof. The
proof occupies the rest of the paper.

\section{Primitive prime divisors} Let $M=(M_n)$ denote the Mersenne
sequence, whose $n$th term is $M_n=2^n-1$. No proof is known
that~$M$ contains infinitely many prime terms. The concept of a
primitive divisor was introduced as a way of showing that new
primes are produced by the terms of $M$, but in a less restrictive
sense. In 1886 Bang~\cite{bang} showed that if~$a$ is any fixed
integer with $a>1$ then the sequence with $n$th term~$a^n-1$ has a
primitive divisor for any index $n>6$. This is a sharp result
because the term~$M_6=63=M_2^2.M_3$ does not have a primitive
divisor. Bang's theorem is remarkable because the number 6 is
uniform across all $a$ and it is small. Indeed, it is not hard to
show that $a=2$ is the only example realizing this bound. Bang's
Theorem was incredibly influential.

In 1892 Zsigmondy~\cite{Z} obtained the generalization that for
any relatively prime $a$ and $b$ with $a>b>0$, the terms of the
sequence $A_n=a^n-b^n$ all have primitive divisors if $n>6$. This
lovely result was re-discovered several times in the early 20th
century and it has turned out to be quite applicable. See~\cite{cp}
and the references therein where applications to Group
Theory are discussed. For example, for fixed~$q$ a prime power,
let $\mathbb F_q$ denote the finite field with~$q$ elements.
Zsigmondy's Theorem applied to the explicit formula for the order
of the group $GL_n(\mathbb F_q)$ shows this order has a primitive
divisor for all large~$n$. Thus Sylow's Theorem can be invoked to
deduce information about the structure of the group.

\begin{definition}\label{def:z} Let~$A=\left(A_n\right)_{n\ge1}$
be an integer sequence. Define
$$
Z(A)= \max\{n\mid A_n\mbox{ does not have a primitive divisor}\}
$$
if this set is finite, and~$Z(A)=\infty$ if not. The number~$Z(A)$
will be called the \emph{Zsigmondy bound} for~$A$.
\end{definition}

Thus Bang's Theorem may be stated: $Z(M)=6$, while Zsigmondy's
Theorem may be stated: $Z(A)\le 6$, where $A_n=a^n-b^n$ as above.
Following Zsigmondy's Theorem, the next major theoretical advance
was made by Carmichael. Let~$u$ and~$v$ denote conjugate quadratic
integers. Consider the integer Lucas sequence~$U$ defined by
$$
U_n=(u^n-v^n)/(u-v).
$$
For example, the Fibonacci sequence $F=(F_n)$ is a Lucas sequence.
Carmichael~\cite{carmichael} showed that if $u$ and $v$ are real
then $Z(U)\le 12$. This too is a sharp result because $F_{12}$
does not have a primitive divisor. The general case was settled by
Bilu, Hanrot and Voutier~\cite{MR2002j:11027}. They proved
that~$Z(U)\le 30$ using a powerful cocktail of methods including
start of the art bounds from Diophantine analysis as well as
massive computations to deal with special cases. Again this is a
sharp result as the sequence generated by the roots of the
polynomial $x^2-x+2$ illustrates. The paper~\cite{MR2002j:11027}
gives details about the long journey from Carmichael's result to
the general case.

\subsection{Elliptic Curves}
Now let $E$ denote an elliptic curve in Weierstrass
form,
$$y^2+a_1xy+a_3y=x^3+a_2x^2+a_4x+a_6
$$
with $a_1,\dots ,a_6 \in \mathbb Z$. The shape of the defining
equation forces the denominator of $x(Q)$ to be an integer square,
for any $Q\in E(\mathbb Q)$. Let $Q\in E(\mathbb Q)$ denote a
non-torsion point. For every $n\in \mathbb N$ write
$$x(nQ)=A_n/B_n^2$$
in lowest terms. The sequence $B=(B_n)$ is an elliptic
divisibility sequence associated to $Q$ and $E$.
Silverman~\cite{silabc} obtained a primitive divisor theorem for
elliptic divisibility sequences arising from curves in Weierstrass
form. It seems likely that a uniform version of this theorem holds
for curves in global minimal form. In other words, if $B=(B_n)$
arises from a rational point on a curve in global minimal form
then~$Z(B)\le N_0$, where $N_0$ is independent of $E$ and $Q$. The
proof of Silverman's Theorem suggests that the Zsigmondy bound is
higher for sequences generated by rational points with small
global canonical height (therefore it is significant that all such
heights are uniformly bounded away from zero). The following
example appeared in~\cite{pi}.
\begin{example}\label{sup}Let $B$ denote the elliptic divisibility sequence generated by the
point $Q=(7107,-602054)$ on the elliptic curve
$$y^2+xy+y=x^3+x^2-125615x+61201397$$ Computations suggest that $Z(B)=39$
and no higher value of $Z(B)$
is known for an elliptic divisibility sequence coming from a
Weierstrass curve in minimal form.
\end{example}
The curve in Example~\ref{sup} was taken from a list of small
height points maintained by Noam Elkies \cite{elkies}. The curves
in Elkies' table are not generally in minimal form but the curve
in Example~\ref{sup} has been rendered in minimal form in order to
estimate~$Z(B)$.

Computations with congruent number curves suggest the Zsigmondy
bound is generally very small. A uniform Zsigmondy bound appears
in~\cite{emw} and~\cite{pi} for an infinite class of sequences
arising from congruent number curves. Specifically, let $T\ge 5$
denote a square-free integer and let $Q$ denote a non-torsion
rational point on the curve
$$y^2=x^3-T^2x.
$$
In~\cite{emw} it was shown that, if $x(Q)<0$ or $x(Q)$ is a
rational square, then $Z(B)\le 21$. In~\cite{pi} this result was
improved by reducing the Zsigmondy bound to $2$ (and allowing any
of $x(Q)$ or $x(Q)\pm T$ to be a square), a bound witnessed by an
infinite family of sequences.
Provided the rank of the curve is positive, there will always be
points satisfying the hypotheses stated: this is because, if
$x(Q)>0$ then $x(Q+[0,0])<0$ and, for any rational point~$Q$,
$$x(2Q) \mbox{ and } x(2Q)\pm T
$$
are all rational squares.

\subsection{Comparisons with the classical theory}

There are notable similarities between the results for elliptic
divisibility sequences and the classical results described
earlier. Both give a uniform bound across infinitely many
sequences which is best possible, both rely upon lower bound on
heights and both reduce the problem to solving a finite number
of Thue-Mahler equations. Also, the uniformity result for Lucas
sequences relies on good bounds from transcendence theory together
with the fact that the answer to Lehmer's problem is known for
quadratic integers. A uniform result for elliptic divisibility
sequences in general would appear to require better elliptic
transcendence results than are currently known,  together with a
proof of Lang's Conjecture on lower bounds for heights of points
on elliptic curves.

In one respect, however, the arithmetic of these two classes of
sequences differs markedly. Bang's Theorem may well have arisen as
part of an attempt to prove the Mersenne Prime Conjecture (which
remains open). On the other hand the analogue of that conjecture
is false for elliptic divisibility sequences on curves of the form
\eqref{homogform}, and heuristics, as well as proofs in special
cases, indicate that it fails for those on Weierstrass curves too,
see~\cite{primeds} and~\cite{pe}.

In the section that follows we will reduce the given problem to
one on the Weierstrass model of a curve birationally equivalent to
the curve~(\ref{homogform}). The method proceeds in a pincer
movement, somewhat similar to that in the two papers~\cite{emw}
and~\cite{pi}. These papers used a good lower bound for the
canonical height of a rational point which were obtained
in~\cite{bst}. Here, our workhorse is the paper~\cite{jed},
although the height bounds are not stated or used in the same way
as in~\cite{jed}. However, beside the similarities, there are many
intriguing differences. Most notably, in this paper we make heavy
use of a numerator sequence on a Mordell curve, see
(\ref{mordellcurve}), for which we can prove a uniform Zsigmondy
bound - see Theorem~\ref{numseq}. Remarkably, this sequence is not
a divisibility sequence, one of the few known cases where a
primitive divisor theorem can be proved for a sequence which lacks
the divisibility property. Also, remarkably, we have been unable
to prove a uniform Zsigmondy bound for the corresponding
denominator sequence on the Mordell curve.

\subsection{The structure of the proof}\label{outline} The proof of
Theorem~\ref{theanswer} relies upon upon two different techniques.
The first one bounds $Z(W)$ above by showing that the
non-existence of a primitive divisor of~$W_n$ implies a certain
divisibility relation (see~(\ref{divstateforV})) involving the
terms~$A_n$, which are defined in~(\ref{def:AnonMordell}). This
relation leads to an inequality which is violated for all
sufficiently large~$n$. This part results in an upper bound
for~$n$ of $n\le 14$.

The second step shows directly that for each of the
indices $n\ge 2$ not covered by the first part, $W_n$ does have a
primitive divisor, by reducing the checking to a finite number of
Thue-Mahler equations. In the reduction, the elliptic division
polynomials, which are elliptic analogues of the cyclotomic
polynomials, play a starring role.

The combination of techniques
described here runs exactly parallel to those used in earlier
primitive divisor theorems such as those of~\cite{bang,
MR2002j:11027, carmichael, Z}.
Note that it is essential to reduce the bound for $Z(W)$ in the
first step as low as possible, in order to keep to a
minimum the number of Thue-Mahler equations which need to be solved
in the second step: without adequate care, computationally infeasible problems result.

The proof uses the well-known bi-rational equivalence
of~(\ref{affinehomogform}) with the Mordell curve
\begin{equation}\label{mordellcurve}
E: Y^2=X^3-432m^2.
\end{equation}
The map is given by
$$u=\frac{36m+Y}{6X}\mbox{ and } v=\frac{36m-Y}{6X}.
$$
If $R\in C(\mathbb Q)$ corresponds to~$Q\in E(\mathbb Q)$ under
the transformation, and we write
\begin{equation}\label{def:AnonMordell}
nQ=\left(\frac{A_n}{B_n^2},\frac{C_n}{B_n^3}\right),
\end{equation}
then
\begin{equation}\label{vnform}\frac{U_n}{W_n}=\frac{36mB_n^3+C_n}{6A_nB_n} \mbox{ and }
\frac{V_n}{W_n}=\frac{36mB_n^3-C_n}{6A_nB_n}.
\end{equation}

The proof of Theorem~\ref{theanswer} exploits both the denominator
sequence $B\nobreak =\nobreak (B_n)$ in~(\ref{def:AnonMordell}),
as in~\cite{emw, pi}, and the numerator sequence $A\nobreak
=\nobreak (A_n)$. The latter may have independent interest.

\begin{theorem}\label{numseq}Let $E$ denote the Mordell curve~(\ref{mordellcurve}) and
suppose $Q$ is a non-torsion point in $E(\mathbb Q)$. Let
$A\nobreak = \nobreak (A_n)$ denote the sequence as defined
in~(\ref{def:AnonMordell}). Then $Z(A) \le 12$.
\end{theorem}

In principle, techniques analogous to those used in~\cite{pi} and
Section~\ref{formssection} can be used here to reduce the bound
stated in Theorem~\ref{numseq}. In practice, however, the
computations involved for some cases are beyond our current
capabilities. We are uncertain about the supremum of the
values~$Z(A)$ as~$m$ varies. Perhaps it occurs when~$m=7$: in this
case~$Z(A)=2$ because~$A_2$ is a proper divisor of~$A_1$, whereas
all the terms~$A_3,A_4,\dots$ have a primitive divisor.

\section{Local arithmetic}\label{local}

Since we are interested in the prime divisors of the numerators
and denominators of $x(nQ)$, for $Q\in E(\Q)$, it makes sense to
address the local arithmetic of $E$.

\textbf{Standing assumption: throughout Section~\ref{local}, $p$
will be a prime other than 2 or 3.}

We will analyze $E(\Q_p)$, as in~\cite[Chapter VII]{aec}, through
the exact sequence
\begin{equation}\label{exact}0\rightarrow E_1(\Q_p)\rightarrow E_0(\Q_p)\rightarrow
\tilde{E}_\textup{ns}(\F_p)\rightarrow 0,\end{equation} where
$E_0(\Q_p)$ is the subgroup of $E(\Q_p)$ consisting of points with
nonsingular reduction modulo $p$, $E_1(\Q_p)$ is the kernel of
reduction modulo $p$, and $\tilde{E}_\textup{ns}(\F_p)$ is the
group of nonsingular points on the curve reduced modulo $p$.
Throughout the paper, the terms {\it nonsingular reduction} and {\it
good reduction} will be used synonymously, as will the terms {\it
singular reduction} and {\it bad reduction}.

For an arbitrary integer sequence $A=(A_n)$, follow Morgan Ward's
terminology \cite{ward} and define the \emph{rank of apparition}
of the prime $p$ in the sequence $A$ to be the least index $n$
such that $p\mid A_n$ (the rank is $\infty$ if no such index
exists). For the sequence $B$ defined above, note that $p\mid B_n$
is equivalent to $\ord_p(x(nQ))<0$, in other words, $nQ\in
E_1(\Q_p)$. The rank of apparition of $p$ in the sequence $B$ is,
then, the order of $Q$ in $E(\Q_p)/E_1(\Q_p)$, and $p\mid B_k$
precisely when $k$ is divisible by this order.  In fact,  the
power to which $p$ divides $B_{kn}$ is entirely predictable once
we know the power to which $p$ divides $B_n$. The following lemma
is obtained in~\cite{aec} by appeal to formal groups (and
by an examination of division polynomials in~\cite{pi}).

\begin{lemma}\label{denompows}
Let $Q\in E_1(\Q_p)$. Then
$$\ord_p(x(kQ))=\ord_p(x(Q))-2\ord_p(k).$$
\end{lemma}

Remember that $p>3$ is a standing assumption in this section.
Lemma~\ref{denompows} is not generally true when $p=2$ (although
it is true when $p=3$).

Note that the commentary above about the rank of apparition leads
to an entirely algebraic interpretation of primitive divisors in
the sequence $B$.  To say that the term $B_n$ fails to have a
primitive divisor is to say that there is no prime $p$ such that
$Q$ has order exactly $n$ in $E(\Q_p)/E_1(\Q_p)$. This quotient,
for primes $p$ of good reduction, is simply $\tilde{E}(\F_p)$. The
rank of apparition of $p$ in the sequence $A$ may similarly be
interpreted in terms of the local arithmetic of $E$, although the
interpretation depends on whether $E$ has good or bad reduction at
$p$.

If $p$ is a prime of bad reduction for $E$, then $p\mid m$
(recalling again that $p>3$ here). The curve $\tilde{E}$ has a
singularity at the point $[0,0]$, and so $p\mid A_n$ precisely
when $nQ$ is singular modulo $p$, in other words, $nQ$ has
non-trivial image in the quotient $E(\Q_p)/E_0(\Q_p)$. Since the
discriminant of $E$ is $\Delta(E)=-2^{12}3^9m^4$, we see that
$\ord_p(\Delta(E))$ is divisible by 4 for all $p\mid m$ other than
3, and one can easily check (see~\cite[Table 15.1]{aec} or use the
addition formula) that
$$E(\Q_p)/E_0(\Q_p)\cong \Z/3\Z.$$
Thus if $p$ appears in the sequence $A$ at all, then it appears in precisely the
terms $A_k$ for which $3\nmid k$.  The primes so occurring are, of course, exactly
the primes of bad reduction dividing $A_1$.

If, on the other hand, $p\nmid m$, then $E(\Q_p)=E_0(\Q_p)$ and
$\tilde{E}=\tilde{E}_\textup{ns}$. Let $H_p\subseteq
\tilde{E}(\bar{\F_p})$ be the subgroup generated by the two points
in $\tilde{E}(\bar{\F_p})$ with $x$-coordinate 0.  It is easy to
verify that $H_p\cong\Z/3\Z$, and clearly $p\mid A_n$ if and only
if  $$n\tilde{Q}\in H_p\setminus\{\Ocal\}$$ where $\tilde{Q}$ is
the image of $Q$ in $\tilde{E}(\F_p)$.  Thus the rank of
apparition of $p$ in $A$ is the order of $\tilde{Q}$ relative to
$H_p\setminus\{\Ocal\}$, or the least $n$ such that $nQ\in
H_p\setminus\{\Ocal\}$.  Note that as $-432m^2= -3(12m)^2$, we
have $$H_p\subseteq\tilde{E}(\F_p)\quad\text{ if and only if
}\quad p\equiv 1\MOD{3},$$ that is, precisely if $E$ has  ordinary
reduction at $p$.  Super-singular primes cannot appear at all in
the sequence $A$, in marked contrast with the situation for the
sequence $B$, in which every prime eventually occurs.  Of course,
even $p$ being a prime of ordinary reduction for $E$ does not
ensure that $p$ has finite rank of apparition in $A$.

The usefulness of Lemma~\ref{denompows}, from our perspective, is
that it allows one to obtain a strong bound on the size of a term
$B_n$ failing to have a primitive divisor, as in~\cite{emw, pi,
is}. The main goal of this section is to prove a similar result
for the sequence $A$.

\begin{lemma}\label{numerpows}
Suppose $\ord_p(x(Q))>0$.  Then, for any $k$ coprime to $3$,
$$\ord_p(x(kQ))=\ord_p(x(Q))+\ord_p(k).$$
\end{lemma}

\begin{proof}
Much of what we need to prove our result has already been
established. If $p$ is a prime of good reduction for $E$, then the
condition $\ord_p(x(Q))>0$ ensures that $\tilde{Q}\in
H_p\setminus\{\Ocal\}$.  In particular, $3\tilde{Q}=\Ocal$, and so
$3Q\in E_1(\Q_p)$.  As $3\nmid k$, we have $k\tilde{Q}\in
H_p\setminus\{\Ocal\}$ as well, and consequently $3kQ\in
E_1(\Q_p)$. Indeed, triplication on $E$ follows the law
\begin{equation}\label{triple}x(3Q)=
\frac{x^9(Q)+2^93^4x^6(Q)m^2+2^{12}3^7x^3(Q)m^4-2^{18}3^9m^6}{9x^2(Q)(x^3(Q)-2^63^3m^2)^2},\end{equation}
and so we deduce (for primes $p\nmid m$) that
\begin{equation*}\ord_p(x(Q))>0\Longrightarrow
\ord_p(x(3Q))=-2\ord_p(x(Q)).\end{equation*}
Applying Lemma~\ref{denompows}, we obtain
\begin{eqnarray*}
2\ord_p(x(kQ))&=&-\ord_p(x(3kQ))\\
&=&-\ord_p(x(3Q))+2\ord_p(k)\\
&=&2\ord_p(x(Q))+2\ord_p(k).
\end{eqnarray*}

The proof is somewhat more involved for the case where $p$ is a prime of
bad reduction.
In this case we know that $E(\Q_p)/E_0(\Q_p)\cong \Z/3\Z$,
and $\ord_p(x(Q))>0$ ensures
that $Q$ is non-trivial in this quotient.  If we let
$$E_n(\Q_p)=\{Q\in E(\Q_p):\ord_p(x(Q))\leq -2n\},$$
then for each $n\geq 0$,
\begin{equation}\label{filtration}E_n(\Q_p)/E_{n+1}(\Q_p)\cong
\Z/p\Z.\end{equation}  This agrees with our definitions above of
$E_0$ and $E_1$, and the fact that $E$ has additive reduction
modulo $p$.  It is worth mentioning that~\eqref{filtration} is
essentially equivalent to Lemma~\ref{denompows}
(although~\eqref{filtration} only holds in the general case for
$n\geq 1$).

Returning to~\eqref{triple} we see that $\ord_p(x(Q))>0$ implies
\begin{equation}\label{badp}
2\ord_p(x(Q))=-\ord_p(x(3Q))+2\ord_p(m).
\end{equation}
If $3Q\in E_n(\Q_p)\setminus E_{n+1}(\Q_p)$ and $r=\ord_p(k)$,
then~\eqref{filtration} tells us that $3kQ\in E_{n+r}(\Q_p)\setminus
E_{n+r+1}(\Q_p)$ so, by~\eqref{badp},
\begin{eqnarray*}
2\ord_p(x(kQ))&=&-\ord_p(x(3kQ))+2\ord_p(m)\\
&=&-\ord_p(x(3Q))+2\ord_p(k)+2\ord_p(m)\\
&=&2\ord_p(x(Q))+2\ord_p(k).
\end{eqnarray*}
%
\end{proof}

\section{The proofs that $Z(A)\le 12$ and $Z(W)\le 14$}\label{mainproof}

After gathering some preliminary results, we proceed with the
proof of Theorem~\ref{numseq}. The proof that $Z(W)\leq 14$, which
advances along similar lines, follows. The step from $Z(W)\leq 14$
to $Z(W)\leq 1$ is taken in Sections~\ref{formssection} and~\ref{threecases}.

\subsection{Preliminaries}
It will not generally be true that the defining equation
(\ref{mordellcurve}) for $E$ is in global minimal form. The
following comes from~\cite[Lemma 1]{jed}.

\begin{lemma}\label{minmodel}If $9\mid m$ write $M=m/9$. Then the
global minimal form $E^*$ for $E$ is
\begin{alignat}{3}
{\rm (I)} \quad &E^*:&\quad y^2 &= x^3-\frac{27m^2}{4} &\text{ if }&2\mid m \text{ and } 9\nmid m\nonumber\\
{\rm (II)}\quad &E^*:&\quad y^2+y &= x^3-\frac{27m^2+1}{4} \quad &\text{ if }& 2\nmid m \text{ and } 9\nmid m\nonumber\\
{\rm (III)}\quad &E^*:&\quad y^2 &= x^3-\frac{3M^2}{4} &\text{ if } &2\mid m \text{ and } 9\mid m\nonumber\\
{\rm (IV)}\quad &E^*:&\quad y^2+y &=x^3-\frac{3M^2+1}{4} &\text{ if
}&2\nmid m \text{ and } 9\mid m.\nonumber
\end{alignat}
The following explicit transformations render the curve~$E$ in
minimal form:
\begin{equation}\label{changecurve}X=u^2x \mbox{ and } Y=u^3y+t
\nonumber\end{equation} where $[u,t]=[2,0], [2,4], [6,0], [6,108]$
(respectively).
\end{lemma}

In the sequel, these four possibilities will be referred to as
Cases (I)-(IV).

\begin{lemma}\label{naivecanon}Write $Q$ for a non-torsion point on $E(\mathbb Q)$, corresponding to
$Q^*\in E^*(\mathbb Q)$, where $E^*$ denotes the minimal model.
Write, for all $n\ge 1$, $x(nQ^*)=a_n/b_n^2$ and $h=\hat h(Q)$.
Also, write~$M=m/9$ if $9\mid m$ and $M=m$ otherwise. Then
$$-\frac23\log M -\frac32\log 3\le hn^2 -\frac18\log
\left|a_n^4+\frac{54M^2a_nb_n^6}{m}\right|\le \frac{1}{12}\log 3.
$$
\end{lemma}

The proof follows immediately from~\cite[Proposition 2]{jed}. Note
the misprint in~\cite{jed} (which has + signs on the left hand
side). We are going to use this in the following form:
\begin{equation}\label{lbforlogan}
hn^2 -\frac{1}{12}\log 3 - \frac18\log
\left|1+\frac{54M^2}{mx_n^3}\right| \le \frac12 \log a_n \mbox{
and }
\end{equation}
\begin{equation}\label{ubforlogan}\frac12\log a_n \le hn^2 +\frac23\log
M +\frac{3}{2}\log 3 .
\end{equation}

\begin{lemma}\label{generallowerbound}Let $P$ denote any non-torsion point
in $E(\mathbb Q)$. Then
\begin{equation}\label{nice}
\hat h(P) \ge \frac{1}{27}\log m -\frac{1}{27}\log 2 -
\frac{1}{36}\log 3> \frac{1}{27}\log m - .0562
\end{equation}
unless $m\equiv \pm 2\MOD 9$ and $m$ has a prime factor
congruent to $1\MOD 6$, in which case
\begin{equation}\label{notnice}
\hat h(P) \ge \frac{1}{27}\log m -\frac{1}{27}\log 2 -
\frac{1}{12}\log 3\ge \frac{1}{27}\log m - .1173
\end{equation}
\end{lemma}

\begin{remark}The difference between the bounds in (\ref{nice})
and (\ref{notnice}) might seem so slight as to be hardly worth
mentioning. However the sieving allowed by the second bound
greatly reduces the amount of manual checking in the sequel.
Having said that, it will become clear in the following proof that
further sieving is possible. The lemma as stated represents a
compromise between further savings on the checking of values
of~$m$ as against a more complicated version of
Lemma~\ref{generallowerbound}.
\end{remark}

\begin{proof}The proof of Lemma~\ref{generallowerbound} uses the analysis
in~\cite{jed}. It begins by estimating a lower bound for $\hat
h(kP)$, for $k=2$ or $3$, then uses $\hat h(kP)=k^2\hat h(P)$ to obtain the bound
sought. The global height is bounded by estimating the local
canonical height $\lambda_p$ at each place. Write
$$\hat h(P)=\sum_{p\le \infty}\lambda_p(P).
$$
We write $Q=kP$ and estimate $\lambda_p(Q)$ for each finite prime $p$.
Suppose that $Q\in E_0(\mathbb Q_p)$, the
non-singular part of the $p$-adic curve; then from~\cite[(2.5)]{jed} we get
$$\lambda_p(Q)=\frac12\log \max \{1,|x(Q)|_p\}-
\frac{1}{12}\log |\Delta^*|_p,
$$
where $\Delta^*$ denotes the discriminant of the minimal equation, so
\begin{equation}\label{nonarch}
\lambda_p(Q) \ge -\frac{1}{12}\log |\Delta^*|_p.\end{equation}

For the archimedean valuation, assume first that $9\nmid m$.
In~\cite[(2.3)]{jed} the following bound is proved:
\begin{equation}\label{arch}\lambda_{\infty}(Q) \ge \frac18\log |x(Q)^4+54m^2x(Q)|-
\frac{1}{12}\log \Delta^*  .\end{equation}
The bound in~(\ref{arch}) holds for any rational point.

If $Q\in E_0(\mathbb Q_p)$ for all primes $p$ then
sum over all $p$, using~(\ref{arch}) and~(\ref{nonarch}), to
obtain
\begin{equation}\label{lowerboundcase1}
\hat h(Q) \ge \frac18\log |x(Q)^4+54m^2x(Q)|=\frac18\log
|x(Q)|+\frac18\log |x(Q)^3+54m^2|,\nonumber
\end{equation}
using the product formula to write $\sum_{p\le \infty}\log
|\Delta^*|_p=0$. Since $$x(Q)^3\ge 27m^2/4,$$ it follows that
$$\hat h(Q) \ge \frac{1}{24}\log \left(\frac{27}{4}\right)+\frac{1}{24}\log m^2 +
\frac18\log m^2 + \frac18 \log \left(\frac{27}{4}+54\right) >
\frac13\log m.
$$
If $m\not\equiv \pm 2\MOD 9$ then, according to~\cite[Page
180]{jed}, $3P\in E_0(\mathbb Q_p)$ for all primes $p$ so we may
put $Q=3P$. The lower bound
$$\hat h(P) > \frac{1}{27}\log m
$$
results, which is stronger than the lower bound~(\ref{nice}). If
$m\equiv\pm 2\MOD 9$ but $p$ has no prime factors congruent to
$1\MOD 6$ then, by ~\cite[Page 180]{jed} again, we may take
$Q=2P$. This time we obtain the lower bound
$$\hat h(P) > \frac{1}{12}\log m
$$
which, again, is stronger than the lower bound in (\ref{nice}).

When $9\mid m$, we use~\cite[(2.4)]{jed}, which gives
\begin{equation}\label{arch9|m}\lambda_{\infty}(Q) \ge \frac18\log |x(Q)^4+6M^2x(Q)|-
\frac{1}{12}\log \Delta^*,\nonumber\end{equation} in place of
(\ref{arch}). Assuming again that $Q\in E_0(\mathbb Q_p)$ for all primes $p$ and summing over $p$, we obtain
\begin{equation}\label{lowerboundcase2}\nonumber
\hat h(Q) \ge \frac18\log |x(Q)^4+6M^2x(Q)|=\frac18 \log
|x(Q)|+\frac18\log |x(Q)^3+6M^2|.
\end{equation}
When $9\mid m$, $x(Q)^3\ge 3M^2/4$ so
$$\hat h(Q) \ge \frac{1}{24}\log \left(\frac{3}{4}\right)+\frac{1}{24}\log M^2 +
\frac18\log M^2 + \frac18 \log \left(\frac{3}{4}+6\right),
$$
which simplifies, upon setting $M=m/9$ to
$$\hat h(Q) \ge \frac13\log m - \frac14\log 3 - \frac13\log 2.
$$
Again by~\cite[Page 180]{jed}, we have $3P\in E_0(\mathbb Q_p)$ for all primes $p$ and putting $Q=3P$ yields
$$\hat h(P) \ge \frac{1}{27}\log m - \frac{1}{36}\log 3 - \frac{1}{27}\log 2 >\frac{1}{27}\log m-.0562
$$
Thus~(\ref{nice}) holds when $9\mid m$.

Finally consider the case when $m\equiv \pm 2\MOD 9$ but $p$ has prime factors congruent to $1\MOD 6$.
In this case $Q=3P\in E_0(\mathbb Q_p)$ for all primes $p\ne 3$.
When $p=3$, use~\cite[(2.6)]{jed} to obtain,
\begin{equation}\label{nonarch3} \lambda_3(Q) \ge -\frac32\log 3-\frac{1}{12}\log |\Delta^*|_3.
\end{equation} Now sum over all $p$, using~(\ref{nonarch}),~(\ref{arch}) and
(\ref{nonarch3}),
to obtain
\begin{equation}\label{lowerboundwithstuff}
\hat h(Q) \ge \frac18\log |x(Q)^4+54m^2x(Q)|-\frac32\log 3.
\end{equation}
Since $x(Q)^3\ge 27m^2/4$,
$$\hat h(Q) \ge \frac{1}{24}\log \left(\frac{27}{4}\right)+\frac{1}{24}\log m^2 +
\frac18\log m^2 + \frac18 \log
\left|\frac{27}{4}+54\right|-\frac32\log 3,
$$
which simplifies to
$$
\hat h(Q) \ge \frac13\log m -\frac13\log 2 - \frac34\log 3.
$$
Writing $Q=3P$ and dividing through by $9$ gives the bound in
(\ref{notnice}).
\end{proof}

\subsection{Proof of Theorem~\ref{numseq}}\label{hereitcomes}

Assume that $n$ is an index such that $A_n$ has no primitive
divisor; the following proof shows that~$n\le 12$. The proof comes
in two steps. For the first step we will show that the assumption
implies a divisibility statement of the following kind:
\begin{equation}\label{divstate}A_n\mid 2^{\mu}3^{\lambda}\rho(n)\prod_{q\mid n}A_{\frac{n}{q}},
\end{equation}
where the product is taken over primes, and
where $\rho(n)$ denotes the product of the primes~$q>3$ which
divide~$n$,
$$\rho(n)=\prod_{3<q\mid n}q.$$ This step is crucial in the
sequel. A little parsimony here greatly reduces the number of
values~$m$ which need to be checked manually.

To prove~(\ref{divstate}), let $p$ be any prime dividing
$A_n$. Assuming~$A_n$ has no primitive divisor, there is a term
$A_k$ with $k<n$ such that $p\mid A_k$.  Initially, suppose that $p\nmid 6$.
If $\alpha$ is the rank of apparition of $p$ in the sequence $A$ then we know,
by the discussion in Section~\ref{local}, that $p\mid A_k$ if and only if
$k=d\alpha$ for some $d$ prime to 3.  Since $\alpha<n$,
there is some prime $q$, necessarily distinct
from $3$, such that $q\alpha\mid n$, and hence $p\mid A_{n/q}$.   Applying
Lemma~\ref{numerpows}, we have
 $$\ord_p(A_n)= \ord_p(A_{\frac{n}{q}})+\ord_p(q)\leq \ord_p(A_{\frac{n}{q}})+1.$$

Now consider the possibilities when $p\mid 6$, subdividing according
to the cases in Lemma~\ref{minmodel}, beginning with the simplest.

\textbf{Case II $\bs2\bs\nmid\bs m, \bs9\bs\nmid\bs m$.}

For the case when $p=2$ notice that there are no rational points
on the minimal model with $x\equiv 0$ mod~$2$. This is because the
expression $(27m^2+1)/4$ is odd and the equation $y^2+y=1$ has no
solutions in $\mathbb F_2$. Hence the maximum value of
$\ord_2(A_n)$ is $2$, using the transformation $X=4x$, and indeed
this is the $2$-adic order of all the~$A_n$. In other words we may
may take $\mu=0$ in~(\ref{divstate}), unless~$n$ is prime, in
which case~$\mu=2$.

Similarly, on the minimal model there are no points with $x\equiv
0\MOD 3$. To see this, notice that $(27m^2+1)/4 \equiv -20\MOD 81$,
if $3\mid m$, and the equation $y^2+y+20\equiv 0\MOD 81$ has
no solutions. If $3\nmid m$ then $(27m^2+1)/4 \equiv 7\MOD 81$
but $y^2+y+7\equiv 0\MOD 81$ has no solutions. Again, using the
transformation $X=4x$, we see that $3\mid A_n$ cannot hold and we may
take $\lambda =0$ in~(\ref{divstate}).

The possibilities for the remaining cases can be summarized as
follows:

\textbf{Case I $\bs2\bs\mid \bs m, \bs9\bs\nmid\bs m$.}
If $4\nmid m$ then $x\equiv 0\MOD 2$ does not hold while if $4\mid m$ then $x\equiv 0\MOD 4$ does
not hold. Using the transformation $X=4x$, we see that $2\le
\ord_2(A_n) \le 3$ so we take $\mu=1$, unless~$n$ is prime, in
which case~$\mu=3$. For the prime $3$, when $3\nmid m$ we find
$x\equiv 0\MOD 3$ does not hold while if $3\mid m$ we find $x\equiv
0\MOD 9$ does not hold. Hence we may take $\lambda =0$ unless
$n$ is prime, in which case $\lambda=1$.

\textbf{Case III $\bs2\bs\mid \bs m, \bs9\bs\mid \bs m$.}
If $4\nmid m$ then $x\equiv 0\MOD 2$ does not hold while if $4\mid m$
then $x\equiv 0\MOD 4$ does
not hold. Using the transformation $X=36x$, we see that $2\le
\ord_2(A_n) \le 3$ so we take $\mu=1$ unless~$n$ is prime, in
which case $\mu=3$. For the prime $3$, we find $x\equiv 0\MOD 3$
does not hold. Hence we may take $\lambda =0$ unless $n$ is prime,
in which case $\lambda=2$ (from $X=36x$).

\textbf{Case IV $\bs2\bs\nmid\bs m, \bs9\bs\mid \bs m$.}
We find $x\equiv 0\MOD 2$
does not hold and neither does $x\equiv 0\MOD 3$. Using the
transformation $X=36x$ gives $\mu=\lambda=0$ unless $n$ is prime,
in which case $\mu=\lambda=2$.

\medskip

Now that~\eqref{divstate} is established, take logarithms:
$$
\log A_n \le \log \rho(n) + \sum_{q\mid n}\log A_{\frac{n}{q}} + \mu
\log 2 + \lambda \log 3,
$$
where the sum runs over primes~$q$. In each case, going to the
minimal model means we must substitute $A_n=u^2 a_n$. This yields
\begin{equation}\label{longineq}
\log a_n \le \log \rho(n) +\sum_{q\mid n}\log a_{\frac{n}{q}}+2\log u
(\omega(n)-1) +\mu \log 2 + \lambda \log 3,
\end{equation}
where
$\omega(n)$ denotes the number of distinct prime divisors of~$n$.

In order to apply Lemma~(\ref{naivecanon}) we will need the
following bound:
\begin{equation}\label{log9}
\log \left|1+\frac{54M^3}{mx_n^3}\right| \le \log 9 = 2\log 3.
\end{equation}
To prove~(\ref{log9}) note firstly that when $9 \nmid m$, we have
$m=M$. Using either form of the minimal equation in
Lemma~\ref{minmodel}, it follows that $x_n^3\ge 27m^2/4$ and the
bound in~(\ref{log9}) follows at once. When $9\mid m$, $M=m/9$ and
Lemma~\ref{minmodel} gives $x_n^3\ge 3M^2/4$. Now~(\ref{log9})
comes out in exactly the same way.

Multiplying~(\ref{lbforlogan}) by $2$ gives
$$2hn^2-\frac{1}{6}\log 3 -\frac14\log
\left|1+\frac{54M^3}{mx_n^3}\right|\le \log a_n
$$
and inserting~(\ref{log9}) gives
$$
2hn^2-\frac23\log 3 \le \log a_n.
$$
Inserting this into~(\ref{longineq}), we get
\begin{equation}\label{evenlongerineq}
2hn^2-\frac23\log 3 \le \log \rho(n) +\sum_{q\mid n}\log
a_{\frac{n}{q}}+2(\omega(n)-1)\log u +\mu \log 2 + \lambda \log 3,
\end{equation}

On the other hand, multiplying~(\ref{ubforlogan}) by $2$,
replacing~$n$ by~$n/q$, and inserting into~(\ref{evenlongerineq})
gives
$$
2hn^2-\frac23\log 3\le \log \rho(n) + \sum_{q\mid
n}\left(2h\left(\frac{n}{q}\right)^2+\frac{4}{3}\log M+ 3\log 3
\right)$$$$+ 2(\omega(n)-1)\log u+\mu \log 2 + \lambda \log 3.
$$
Re-arranging gives
$$2hn^2\left(1-\sum_{q\mid n}\frac{1}{q^2}\right)\le \log \rho(n)+\frac{4}{3}\omega(n)\log M$$$$
+ 2(\omega(n)-1)\log u +\mu \log 2 + \left(\lambda
+3\omega(n)+\frac23\right)\log 3.
$$
To ease the notation, write
$$f(n)=1-\twolinesum{q\mid n}{q\ {\rm prime}}\frac{1}{q^2}.
$$
Since $M\le m$ in each case, dividing by $\log m$ yields

\begin{equation}\label{generic}
\frac{2hn^2}{\log m}f(n)\le \frac{4}{3}\omega(n)+ \frac{\log
\rho(n)+\omega(n)\log (27u^2) +\log
(2^{\mu}3^{\lambda+\frac23}/u^2)}{\log m}.
\end{equation}

\textbf{Case II $\bs2\bs\nmid\bs m$, $\bs9\bs\nmid\bs m$.}

In this case, $u=2$ and
$\mu=\lambda=0$ unless $n$ is prime, in which case $\mu=2$ and
$\lambda=0$. When~$n$ is composite, assume $m \ge 40$. The bound
(\ref{nice}) inserted into~(\ref{generic}) now forces $n\le 12$.
Assuming $m \ge 290$, the bound~(\ref{notnice}) inserted into
(\ref{generic}) forces $n\le 12$. Thus we need to check all cases
manually when $m\le 40$ and all cases $m\le 290$ when $m\equiv \pm
2\MOD 9$ and $m$ has a prime divisor $\equiv 1$ mod~$6$. This
will follow in the next section.

When~$n$ is prime, $\omega(n)=1$ but we need to take $\mu=3$ and
$\lambda=1$. When $m\ge 40$,~(\ref{nice}) and~(\ref{generic})
force $n\le 7$. When $m \ge 290$,~(\ref{notnice}) and
(\ref{generic}) force $n\le 7$.

\begin{remark}To obtain these explicit upper bounds rigorously requires a
little more than checking~(\ref{generic}), because some of the
functions that appear in~(\ref{generic}) are not monotonic. In
fact~(\ref{generic}) implies a weaker inequality, where
$\omega(n)$ is replaced by its upper bound of $\log n/\log 2$,
$\rho(n)$ is replaced by its upper bound of~$n$, and $f(n)$ is
replaced by its lower bound of $.547$. The resulting inequality
yields a rigorous upper bound and the remaining cases can be
checked manually. In none of the cases needed does this higher
bound exceed $24$ so the amount of extra checking is negligible.
\end{remark}

Now we summarize the arguments in the other three cases:

\textbf{Case I $\bs2\bs \mid \bs m$, $\bs9\bs\nmid\bs m$.} Here $u=2$ and $\mu=1$ and
$\lambda=0$ unless $n$ is prime, in which case $\mu=3, \lambda
=1$. Applying Lemma~\ref{generallowerbound} together with
(\ref{generic}) as before gives the same bounds as the previous
case with the same manual checking to be done.

\textbf{Case III $\bs2\bs \mid \bs m$, $\bs9\bs \mid \bs m$.} Here $u=6$, $\mu=1$ and
$\lambda=0$ unless $n$ is prime, in which case, $\mu=3$ and
$\lambda =2$. Inserting Lemma~\ref{generallowerbound} into
(\ref{generic}), the same bounds arise.

\textbf{Case IV $\bs2\bs \nmid\bs m$, $\bs9\bs \mid \bs m$.} Here $u=6$ and $\mu=\lambda=0$
unless $n$ is prime, in which case $\mu=2, \lambda =2$. Inserting
into~(\ref{generic}) and using Lemma~\ref{generallowerbound} gives
exactly the same bounds as before.

The proof of Theorem~\ref{numseq} is complete subject to checking
various values of~$m$; see Appendix~\ref{smallm} for the details.
\hfill\qed

\subsection{Proof that $\bs Z\bs(\bs W\bs)\bs\le \bs1\bs4$}\label{hereitcomesZW}

Note that if a prime $p$ is a primitive divisor of a term $A_n$,
where $n\geq 2$, then it cannot divide~$B_r$ with $r<n$. We have
seen that a prime of bad reduction which divides $A_n$, with
$n>1$, also divides $A_1$, and so such a prime cannot be primitive
divisor of $A_n$. If $p\geq 5$ is a primitive divisor of $A_n$
then $p$ is a prime of good reduction for $E$ and so, by the
results in Section~\ref{local}, the rank of apparition of $p$ in
the sequence $B$ is $3n$. When $p\mid 6$, the same conclusion
follows by manual checking, using the same case-by-case analysis
as in Section~\ref{hereitcomes}.

The proof that $Z(W)\le 14$ runs along almost the same lines as
the proof that~$Z(A)\le 12$. Suppose $p\geq 5$ is a primitive
divisor of $A_n$, and suppose that $p^a\| A_n$. Then either $p$
divides $W_n$, in which case it is a primitive divisor, or
$p^a\mid 36mB_n^3+C_n$, according to~\eqref{vnform}. By the
identity
$$(36mB_n^3+C_n)^3+(36mB_n^3-C_n)^3=m(6A_nB_n)^3,$$
which is simply~\eqref{homogform}, we have $p^a\mid 36mB_n^3-C_n$,
whence $p^a\mid 36mB_n^3$.  As $A_n$ is prime to $B_n$, and as
$p\geq 5$, we have $p^a\mid m$.  Thus $W_n$ fails to have a
primitive divisor just in case the `primitive part' of $A_n$
divides $m$, or
\begin{equation}\label{divstateforV}A_n\mid m2^{\mu}3^{\lambda}\rho(n)\prod_{q\mid n}A_{\frac{n}{q}}.
\end{equation}

The arguments of the previous section now apply {\it mutatis
mutandis} to show that
\begin{equation}\label{genericforV}
\frac{2hn^2}{\log m}f(n)\le 1+\frac{4}{3}\omega(n)+ \frac{\log
\rho(n)+\omega(n)\log (27u^2) +\log
(2^{\mu}3^{\lambda+\frac23}/u^2)}{\log m}.
\end{equation}
When~$n$ is composite, assume $m \ge 40$. The bound~(\ref{nice})
inserted into~(\ref{genericforV}) now forces $n\le 14$ in all
cases. Assuming $m \ge 290$, the bound~(\ref{notnice}) inserted
into~(\ref{genericforV}) forces $n\le 14$. This completes the
proof that $Z(W)\le 14$, subject to checking various values
of~$m$. The values of $m$ required to complete the argument are checked in Appendix~\ref{smallm}.\hfill\qed

Thus the proof of Theorem~\ref{theanswer} follows once we
establish that $W_n$ has a primitive divisor for all $2\le n\le
14$. To that end we now turn.


\section{Primitive divisors in specific terms}\label{formssection}

In the present section we concern ourselves with the question of,
for fixed $n$, which sequences $W$ fail to have a primitive
divisor in the $n$th term.  For many small values of $n$ we can
show that there are no such sequences, and this will bridge the
gap between the work in Section~\ref{mainproof} and the goal of
Theorem~\ref{theanswer}.  To ease notation, we will frequently
write $D$ in place of $-432m^2$.

\begin{proposition}\label{smallnprop}
Let $W$ be a sequence as defined above.  Then $W_n$ has a
primitive divisor for each $2\leq n\leq 14$.
\end{proposition}

The proof will come in several pieces.  For $n\geq 5$ not
divisible by 3, the methods developed in~\cite{pi} suffice to
treat this problem.  Although we make use of the special form of
$D=-432m^2$, much of the argument will work for general Mordell
curves.  Indeed, we reduce the proposition to the checking of
finitely many elliptic divisibility sequences arising from Mordell
curves, none of which turn out to be curves of the special form
under consideration.  When $n\geq 6$ is divisible by 3 some
problems arise, but only a slight modification of the method of~\cite{pi}
is needed.  The cases $n=4$, $n=3$, and $n=2$,  treated
in Section~\ref{threecases}, are dispatched largely through
\emph{ad hoc} means, although the spirit of the proof remains the
same.

\subsection{Division polynomials}

We exploit, as above, the short Weierstrass form of the equation,
and will in fact show that $W_n$ has a primitive divisor coprime
to $6m$ for the $n$ listed in the lemma. We will consider the
curves $C$ and $E$, as defined in~\eqref{homogform} and~\eqref{mordellcurve}
respectively, bi-rationally equivalent by the
map defined in~\eqref{vnform}.  There are (see, for example,~\cite{aec})
rational maps $\phi_{n}, \psi_{n}, \omega_{n}\in
\Q(E)$, the function field of $E$, such that, for each $n\in\Z$
and $Q\in E(\Q)$,
$$nQ=\left(\frac{\phi_{n}(Q)}{\psi_{n}^2(Q)}, \frac{\omega_{n}(Q)}{\psi_{n}^3(Q)}\right)=
\left(\frac{A_n}{B_n^2}, \frac{C_n}{B_n^3}\right).$$
In Section~\ref{mainproof}, we exploited the existence of
primitive divisors in the sequence $A$.  Here we use the
properties of the division polynomial $\psi_n$, much as in~\cite{pi},
to establish the existence of primitive divisors in
$W_n$.  An approach employing the polynomials $\phi_n$ would be
similarly successful, and would be in keeping with the flavour of
Section~\ref{mainproof}, but such an approach also encounters
serious computational difficulties.

We will show that for all $n$ under consideration (save those
treated later as special cases), there is a prime $p$ dividing
$B_n$ such that $p\nobreak \nmid \nobreak 6A_kB_k$ for any $k<n$.
Such a prime will divide $W_n$, as $\gcd(B_n, C_n)=1$, but not
$W_k$ for $k<n$. That is, $p$ will be a primitive divisor of
$W_n$.  It turns out that this amounts to showing that there is
some prime $p\nmid\Delta(E)$ such that $\ord_p(F_n(Q))>0$, for a
certain $F_n\in\Q(E)$ constructed below.  The hypothesis that
there is no such prime will lead us to a non-trivial solution of a
certain Thue-Mahler equation which depends on $n$ but, critically,
not on $Q$ or $m$. By explicitly solving the Thue-Mahler equations
in question, we will treat all bar three cases of
Proposition~\ref{smallnprop}.

As $\Q(E)=\Q(x, y)$, with $y^2=x^3+D$ (for fixed $D=-432m^2$), we may write $\phi_{n}$, $\psi_{n}$,
and $\omega_{n}$ as polynomials in $x$ and $y$, where $y$ occurs only to the first power.
In fact (see~\cite{aec}), we may write $\phi_{n}$ as a polynomial in $x$ alone, and either
$\psi_{n}$ or $y^{-1}\psi_{n}$ similarly, as $n$ is odd or even respectively. Note that, in
all cases, $\psi_{n}^2$  may be written as a polynomial in $x$, and we will view $\phi_n$
and $\psi_n^2$ as elements of $\Q[x]$.  In fact, if $E[n]$ denotes the kernel of multiplication
by $n$ in $E(\bar{\Q})$,
\begin{equation}\label{divispoly}\psi_{n}^2(Q)=n^2\prod_{\substack{T\in E[n]\\ T\neq\Ocal}}
(x(Q)-x(T)).\end{equation} Note that each linear term on the right
side occurs precisely twice, except those corresponding to $T\in
E[2]$.

As our proof relies on the properties of $\psi_n$, it is incumbent
upon us to make a few observations (these remarks are made as well
in~\cite{pi}). Note that if $3\nmid n$, then $(x, y)\in E[n]$
implies both $x\neq 0$ and
$$(x, y), (\zeta x, y), (\zeta^2x, y)\in E[n]$$
where $\zeta$ is a primitive 3rd root of unity.  In particular,
$\psi_n$ or $\psi_n/y$, depending on the parity of $n$, is a
polynomial in $x^3$; $\psi_n^2\in\Q[x^3]\subseteq\Q[x]$. When
$3\mid n$, on the other hand,  $\psi_n^2\in x^2\Q[x^3]$.

Although $\psi_{n}$ depends on $m$, our stated aim is to
construct, from each term $W_n$ without a primitive divisor, a
solution to a Thue-Mahler equation that is independent of $m$.
Note that the points on $y^2=x^3+D$ and those on $y^2=x^3+D'$ are
related by the scaling map
$$\left(X, Y\right)\leftrightarrow\left(X\left(\frac{D'}{D}\right)^\frac{1}{3},
Y\left(\frac{D'}{D}\right)^\frac{1}{2}\right).$$ Thus, the
dependence of~\eqref{divispoly} on $D$ is transparent. Along the
same lines as~\cite{pi}, we observe that $\psi_{n}^2$ may be
written, over $\Z$, as a binary form in $x^3$ and $4D$.  We will
abuse notation somewhat, and denote this form by $\psi_n^2$ as
well, so that
$$\psi_n^2(Q)=\psi_n^2(x^3, 4D).$$

In general, $\psi_n^2\in\Q(E)$ is not irreducible.  Aside from
being a square 
when $n$ is odd,
the function has several obvious
factors. Let $e_\infty(T)$ denote the order of the torsion point
$T\in E(\bar{\Q})$, and let
\begin{equation}\label{eff}F_n^2(Q)=
\epsilon^2(n)\prod_{\substack{T\in
E(\bar{\Q})\\e_\infty(T)=n}}(x(Q)-x(T)),\end{equation}
where $$\epsilon(n)=\begin{cases}p & \text{if }n=p^a\text{ is a prime power}\\
1 & \text{otherwise.}\end{cases}$$
It is clear from~\eqref{divispoly} that
$$\psi_n^2=\prod_{d\mid n}F_d^2.$$
Note that, when $n\neq 2$, each term on the right of~\eqref{eff} occurs precisely twice,
allowing us to define $F_n\in\Q[x]$ implicitly in this way.

The functions $F_n$ may be viewed as the elliptic analogues of the
cyclotomic polynomials exploited in~\cite{MR2002j:11027,
carmichael, schin} to treat the analogous problem for Lucas
sequences. By the same arguments as in the case of $\psi_n$, the
product in~\eqref{eff} defines a binary form over $\Z$ in $x^3$
and $4D$, at least when $n\neq 3$.  We will, then, write $F_n(x^3,
4D)$ for the product above.  We will frequently pass between the
rational function $F_n(Q)$ and the binary form $F_n(A_1^3,
4DB_1^6)$.  It is worth noting that, as $A_1$ and $B_1$ are
co-prime, the primes appearing to a positive power in $F_n(Q)$
occur to the same power in $F_n(A_1^3, 4DB_1^6)$, except possibly
those dividing $\epsilon(n)$.

\subsection{Division polynomials in finite fields}

Note that if $p$ is a prime of good reduction for $E$, then the
equations from~\cite{aec} define the division polynomials in the
same way for $E(\F_p)$.  In particular, if $e_p(Q)$ is the order
of the image of the point $Q\in E(\Q)$ in the group $E(\F_p)$,
then $e_p(Q)$ is precisely the rank of apparition of the prime $p$
in the sequence $(\psi_n(Q))_{n\geq 1}$, that is, the smallest $n$
such that $\ord_p(\psi_n(Q))>0$.  Equivalently, $e_p(Q)$ is the
unique $n$ such that $\ord_p(F_n(Q))>0$.  Note, on the other hand,
that $x(Q)\equiv 0\MOD{p}$ implies that $Q$ is a point of order 3
in $E(\F_p)$.  More generally, then, if $n$ is the rank of
apparition of $p$ in the sequence $(\phi_n(Q))_{n\geq 1}$, it
follows that $e_p(Q)=3n$.  The converse may fail, of course, as
there may be points of order 3 in $\tilde{E}(\F_p)$ other than
those with $x=0$. The following observation will be useful.
\begin{lemma}\label{ayadlem}
Let $Q\in E(\Q)$, let $A_n$ and $B_n$ be defined as above, and let
$p\nmid\Delta(E)$. Then the following hold:
\begin{enumerate}
\item if $\ord_p(\psi_n(Q))>0$, then $p\mid B_n$;\label{parta}
\item  if $\ord_p(\phi_n(Q))>0$, then $p\mid A_n$; \item if
$\ord_p(F_n(Q))>0$ and $3\nmid n$, then $p\mid B_n$ and $p\nmid
6A_kB_k$ for $k<n$.
\end{enumerate}
\end{lemma}

\begin{proof}
Note that results similar to (a) and (b) are derived in~\cite{ayad, pi, is}.
In particular, the result in~\cite{ayad} is
stronger: the hypothesis is merely that the reduction of $Q$
modulo $p$ is not singular.  The lemma is demonstrated here for
completeness.

As
$$\frac{A_n}{B_n^2}=\frac{\phi_n(Q)}{\psi_n^2(Q)},$$
the conclusions of (a) and (b) can only fail if
$\ord_p(\phi_n(Q))$ and $\ord_p(\psi_n(Q))$ are simultaneously
positive.  If this is the case, then $\ord_p(\psi_k(Q))>0$ for
$k=n\pm 1$ as, by definition (see~\cite{aec}),
$$\phi_n=x\psi_n^2-\psi_{n-1}\psi_{n+1}.$$

By an easy modification of Lemma 4.1 of~\cite{ward} (for details
on this modification see~\cite{is}), we have $\ord_p(\psi_3(Q))>0$
and $\ord_p(\psi_4(Q))>0$.  Computing the resultants of the binary
forms $\psi_3(x^3, 4D)$ and $\psi_4^2(x^3, 4D)$ in $\Z[D]$, we see
that this can happen only if $p\mid 6m$.  That is, if $p\mid
\Delta(E)$.

Now suppose $\ord_p(F_n(Q))>0$.  Part (a) ensures that $p\mid
B_n$, and the image of $Q$ in $E(\F_p)$ has order exactly $n$.  If
$p\mid B_k$ for $k<n$, then $\ord_p(\psi_k(Q))>0$, and hence
$e_p(Q)\leq k<n$, a contradiction.  Similarly, if $p\mid A_k$ then
we have $\ord_p(\phi_k(Q))>0$.  Hence the prime $p$ has finite
rank of apparition in the sequence $(\phi_n(Q))_{k\geq 0}$, and so
$n=e_p(Q)$ is divisible by 3. We assumed it was not.
\end{proof}

\subsection{Thue-Mahler equations}

The factors $F_n$ are the binary forms exploited in~\cite{pi, is} to show the existence of
primitive divisors in specific terms of elliptic divisibility sequences arising from elliptic
curves in short Weierstrass form.  We will use a similar approach here to construct primitive
divisors of the sequence $B$.
The following, simple observation will be used repeatedly.
\begin{claim}\label{theclaim}
Let
$$s=A_1^3/\gcd(A_1^3, 4D), \quad t=4DB_1^6/\gcd(A_1^3, 4D).$$
Then for all primes $p\geq 5$ dividing $\gcd(A_1, D)$, we have
$$\ord_p(s)>\ord_p(t)=0.$$
Furthermore, for all primes $p\geq 5$, the quantity $\ord_p(t)$ is even.
\end{claim}

\begin{proof}[Proof of the Claim]
Note that, as $m$ is cube-free, we have $$\ord_p(m)\in\{1, 2\}$$
for all $p\geq 5$ dividing $m$. Suppose that $p\mid A_1$ (and
hence $p\nmid B_1$).  If $\ord_p(m)=1$, then
$$\ord_p(A_1^3)\geq 3 > 2 =\ord_p(4DB_1^6)=\ord_p(-432m^2).$$
Thus $\ord_p(s)>\ord_p(t)$.

If, on the other hand, $\ord_p(m)=2$, then $p^3$ divides the
right-hand-side of $$C_1^2=A_1^3-432m^2B_1^6.$$ But then $p^4$
divides $A_1^3-432m^2B_1^6$, as $p$ divides the left-hand-side of
the above to an even power.  Now $p^2\mid A_1$, hence
$$\ord_p(A_1^3)\geq 6>4=\ord_p(-432m^2),$$
and again $\ord_p(s)>\ord_p(t)$.
The fact that $\ord_p(t)=0$ simply follows from $s$ and $t$ being, by construction, relatively prime.

Now suppose that $p\geq 5$, and further that $\ord_p(t)\neq 0$.  We have just shown that this
ensures $p\nmid A_1$, and so $$\ord_p(t)=\ord_p(-432m^2)=2\ord_p(m)\equiv 0\MOD{2}.$$
\end{proof}

We are now in a position to present the main tool in the proof of
Proposition~\ref{smallnprop} for the values of $n\geq 5$ which are
not divisible by $3$.  The values $n=2, 4$, and those values
divisible by 3, require a slightly more careful treatment.

\begin{lemma}\label{tmlem}
Suppose that $W_n$ has no primitive divisor, for $n\geq 5$ not
divisible by 3.  Then, for $s$ and $t$ defined as above, $F_n(s,
t)$ is divisible only by primes dividing $6\epsilon(n)$.
\end{lemma}

\begin{proof}
We have seen that if $p\nmid \Delta(E)$ and $\ord_p(F_n(Q))>0$,
then $p$ is a primitive divisor of $W_n$.  Thus, if $W_n$ has no
primitive divisor, we must have $\ord_p(F_n(Q))>0$ only for
$p\mid\Delta(E)$.  Clearing denominators, we see that the only
primes dividing the integer $F_n(A_1^3, 4DB_1^6)$ are those of bad
reduction for $E$ or, possibly, those dividing $\epsilon(n)$.  If
$p\nmid6\epsilon(n)$ is one such prime, then $p\mid D$.  But then
$$F_n(A_1^3, 4DB_1^6)\equiv \epsilon(n)A_1^{3\deg{F_n}}\MOD{p}.$$
In particular, we can only have $p\mid F_n(A_1^3, 4DB_1^6)$ if
$p\mid A_1$.  In the latter case we have $p\mid\gcd(A_1, D)$
whence, by Claim~\ref{theclaim}, we obtain
$\ord_p(s)>\ord_p(t)=0.$ Note that $F_n(0, 1)=\pm 1$ as
in~\cite{pi} (see also the tables in Section~\ref{divpolytables}).
Thus
$$F_n(s, t)\equiv\pm t^{\deg{F_n}}\not\equiv 0\MOD{p}.$$
In particular, the prime divisors of $F_n(s,t)$ are at most those of $6\epsilon(n)$:
$$F_n(s,t)=\pm 2^\alpha 3^\beta \epsilon(n)^\gamma.$$
\end{proof}

Note, for the purposes of impending computations, that we can say
somewhat more.  If $n=p^a$ is a prime power (and so
$\epsilon(n)=p$), then $\ord_p(F_n(Q))>0$ and $p\mid\Delta(E)$
together imply $p\mid m$. In this  case $D\equiv 0\MOD{p^2}$, ergo
$$F_n(A_1^3, 4DB_1^6)\equiv pA_1^{3\deg(F_n)}\MOD{p^2}.$$
Thus if $p\nmid A_1$, we have $\ord_p(F_n(A_1^3, 4DB_1^6))=1$. If,
on the other hand, $p\mid A_1$, we have $\ord_p(s)>\ord_p(t)$, and
so $p\nmid F_n(s, t)$.  In the notation of the proof, then, we may
take $\gamma\in\{0, 1\}$.  We will see, by examining the
individual forms, that we can control the exponents of 2 and 3 by
elementary means as well.

\subsection{Solving the Thue-Mahler equations}\label{thuesolve}
We have reduced the proof of Proposition~\ref{smallnprop} to
treating the special cases $n=2, 4$ and $3\mid n$, as well as
solving a number of Thue-Mahler equations. Although the proof
appears to require us to find all solutions to the Thue-Mahler
equations
$$F_n(s, t)=2^\alpha3^\beta\epsilon(n)^\gamma,$$
we have already seen that we need only consider equations wherein
$\gamma\in\{0,1\}$.  Although we
can, \emph{a priori}, restrict the exponents $\alpha$ and $\beta$
by employing such techniques as
lower bounds on linear forms in $p$-adic logarithms, it  turns out
that we may also do so by more elementary means.

For the various $n$ under consideration, consider $F_n\MOD{2}$,
and note that
$$F_n(s, t)\equiv 1\MOD{2}$$
whenever $\gcd(s, t)=1$ (note that these forms are available below for examination).  Similarly, we may
consider $F_n\MOD{3}$ and reduce the possible values of $\beta$ to three: 0, $\deg(F_n)$,
or $\frac{3}{2}\deg(F_n)$.  Thus, our Thue-Mahler equations may be reduced to 6 or 12 Thue equations
(depending on whether or not $n$ is a prime power) of the form
$$F_n(s, t)=(-1)^\delta 3^\beta \epsilon(n)^\gamma,$$
with $\delta, \gamma\in\{0, 1\}$ and $\beta\in\{0, \deg{F_n}, \frac{3}{2}\deg(F_n)\}$.

We may also disregard several possible solutions in advance.  For
example, some of these equations possess a solution $(s, t)$ with
$t\nobreak = \nobreak 0$.  Such a solution cannot arise from a
pair $(s, t)$ as constructed in Claim~\ref{theclaim}, however, as
this would necessitate  either $m=0$ or $B_1=0$. Furthermore, we
have from Claim~\ref{theclaim} that $\ord_p(t)$ is even for all
$p\geq 5$, and so we may ignore solutions $(s, t)$ that fail to
have this property.  As
\begin{gather*}
t\gcd(A_1^3, 4D)=-432m^2B_1^6<0\\
\intertext{and}
(4s+t)\gcd(A_1^3, 4D)=4(A_1^3+DB_1^6)=4C_1^2>0,
\end{gather*}
we may conclude as well that $t<0$ and $4s+t>0$ (here, notice that
$C_1=0$ only for points $Q$ of order 2, which are not the type
under consideration).  Any solution not satisfying these
inequalities may be discarded as well.  Solutions wherein $s=0$
correspond to points $Q$ of order 3, and the solution $(s, t)=(1,
-1)$ gives rise to the point $(12, 36)$, a point of order 3 on
$y^2=x^3-432$. We shall call solutions falling into the above
categories \emph{expected}.

A computation in PARI/GP~\cite{PARI2} shows that there are no
solutions to any of the Thue equations above, other than
(possibly) these expected solutions.  In the appendix, we list the
binary forms $F_n$ for the various values of $n$, in order that
the reader may confirm these findings. We should note that the
Thue equation solver in PARI/GP assumes, by default, the truth of
the Generalised Riemann Hypothesis.  This default was overridden,
and our results verified unconditionally. The most strenuous
computation arises in the case $n=11$, in which a Thue-Mahler
equation of degree 20 must be treated. (The binary form arising in
the case $n=13$ is of greater degree but factors, and so the
equations above may be treated by elementary means.)  Even this
computation, however, took well less than a minute (on a 1.83 GHz
MacBook with 512MB of RAM).

\subsection{The case $3\mid n$, but $n\neq 3$}
Some care must be taken when $3\mid n$ and $n\neq 3$, but the
methods are not fundamentally different.

If $F_n\in\Q(E)$ is defined as above, where $n=3k$, and if $p$ is
a primitive divisor of $A_k$, then we have $\ord_p(F_n(Q))>0$.  It
is not clear, then, that the primitive divisor of $B_n$ which we
prove to exist by the method  above fails to divide $A_k$.  We
must modify our argument to show that $B_n$ has a primitive
divisor distinct from those coming from $A_k$.  The binary form
$F_n$ factors in this case, and we will show that it suffices to
show that one of the factors has a prime divisor other than 2 or
3.

Let $H\subseteq E(\bar{\Q})$ be the group
$$\left\{(0, \sqrt{D}), (0, -\sqrt{D}), \Ocal\right\}\cong\Z/3\Z.$$
Note that for $D=-432m^2$, $H\subseteq E(\Q(\sqrt{-3}))$, although
much of what is written here
applies for all $D$.
For $Q\in E(\bar{\Q})$, let $e(Q, H)$ denote the least $k$ such that
$kQ\in H$, if one exists,
and $e(Q, H)=\infty$ otherwise.  Set, for $n$ divisible by 3,
$$G_n^2(Q)=\prod_{\substack{T\in E(\bar{\Q})\\ e(T, H)=n/3}}(x(Q)-x(T)).$$
As $H$ is fixed by the action of Galois on $\bar{\Q}$, we have
immediately that $G_n\in \Q(E)$. If $G_n$ vanishes at $Q$, then
$A_n=0$, and so we see that $G_n\mid \phi_n$.  If we define
$H_p\subseteq E(\bar{\F_p})$ to be the analogous subgroup (for
$p\nmid\Delta(E)$), and $e_p(Q, H)$ to be the analogous value,
then $e_p(Q, H)$ is the rank of apparition of $p$ in the sequence
$A$, as discussed above.   Arguments similar to those before show
that $G_n$ may be written as a binary form in $x^3$ and $4D$ with
coefficients in $\Z$, and we will denote this form by $G_n(X, Y)$.
For example,
$$G_9(X,Y) = X^3-24X^2Y+3XY^2+Y^3.$$
Note that the roots of $G_n$ are points on $E(\bar{\Q})$ of order precisely $n$.
Thus $G_n$ divides $F_n$.
Let $$F_n(X,Y)=G_n(X,Y)\tilde{F}_n(X,Y).$$

\begin{lemma}
Suppose that $W_n$ has no primitive divisor, for $3\mid n$.  Then
$\tilde{F}_n(s,t)$ is a $\{2, 3\}$-unit.
\end{lemma}

\begin{proof}
Suppose that $p\nmid\Delta(E)$ and $\ord_p(\tilde{F}_n(Q))>0$. It
follows that $\ord_p(F_n(Q))\nobreak >\nobreak 0$, and so, just as
in Lemma~\ref{ayadlem}, we have $p\mid B_n$ and consequently
$p\mid W_n$.  We have, by the same argument as is
Lemma~\ref{ayadlem}, that $p\nmid B_k$ for $k<n$. Suppose $p\mid
A_k$, for some $k<n$, and suppose without loss of generality that
$k$ is the least such value.  Then we have $e_p(Q)=3k$, and so
$n=3k$.  We have as well that $e_p(Q, H)=n/3$, so
$\ord_p(G_n(Q))>0$.  But $\tilde{F_n}$ and $G_n$ can have no
common roots modulo $p$. If $Q$ were such a root, then
$x(n\tilde{Q})\equiv 0\MOD{p}$ while also $x(n\tilde{Q})^3\equiv
-4D\MOD{p}$, clearly impossible if $p\nmid D$.  So we must not
have $p\mid A_k$.

If, on the other hand, $p\mid\Delta(E)$ and $p\geq 5$, we have
$p\mid m$. As above, we have $\ord_p(s)>\ord_p(t)=0$ and so, as
$\tilde{F}_n(1, x)$ is monic, $p\nmid \tilde{F}_n(s,t)$. The
result is proved.
\end{proof}

Note that the lemma is true in the case $n=3$.  The binary form
$\tilde{F}_3$ has degree 1, however, and so we cannot proceed in
the same way as we will for $n\geq 6$. As in the case of
Lemma~\ref{tmlem}, some care must be taken with the exponents.  We
note that $\tilde{F}_n(s, t)\equiv 1\MOD{2}$ for all relatively
prime $s$ and $t$, while
$$\ord_3(\tilde{F}_n(s, t))=\begin{cases}  0, 3,\text{ or } 5 & \text{if }n=6\\
0, 9, \text{ or }13 & \text{if }n=9\\
0, 12, \text{ or } 18 & \text{if }n=12.
\end{cases}$$

Computations in PARI/GP  are as in the previous case, and reveal no unexpected solutions.

\section{Three special cases}\label{threecases}

We have, thus far, shown that $W_n$ has a primitive divisor for
all $n\geq 5$.  We treat the cases $2\leq n\leq 4$ here,
completing the proof of Theorem~\ref{theanswer}.

\subsection{The cases $\bs n\bs=\bs3$ and $\bs n\bs=\bs4$}
To treat the case $n=3$, we show that $A_3$ must have a prime
divisor not dividing $6mA_1A_2$.   Such a divisor, being a prime
of good reduction for $E$, must divide $W_3$, but cannot divide
$W_1$ or $W_2$. Suppose, to the contrary, that $A_3$ has no prime
divisor other than those dividing $6mA_1A_2$. Comparing $\phi_3$
to $\psi_2$ and $\phi_2$, we see that then any prime $p$ with
$\ord_p(\phi_3(Q))>0$ must be a prime of bad reduction for $E$. If
$p\geq 5$ is a divisor of $m$ we see, just as above, that
$$\phi_3(s, t)=s^3-24s^2t+3st^2+t^3$$ is not divisible by $p$.  So
$\phi_3(s,t)$ is a  $\{2, 3\}$-unit.  Exactly as in previous
cases, we reduce this to
$$s^3-24s^2t+3st^2+t^3=\pm 3^\beta,$$
with $\beta\in\{0, 3, 4\}$, by considering the possible values of $\phi_3$ modulo 2 and 3.
Computation produces no unexpected solutions.

For the case $n=4$, we will show that $A_4$ has a prime divisor
not dividing $6mA_1A_2A_3$.  As this divisor is a prime of good
reduction for $E$, it divides $W_4$, witnessing that $W_4$ has a
primitive divisor.  We  consider one particular factor of
$\phi_4^2(X,Y)$. Let
$$F^*_4(X,Y)=X^4-134X^3Y-84X^2Y^2-32XY^3-2Y^4.$$
One can verify that this binary form divides $\phi_4^2$, and
(by the computation of some resultants)
that $F^*_4(A_1^3, 4DB_1^6)$ has no prime factors in common with
$A_1, A_2, A_3, B_1, B_2, B_3$,
save possibly some divisors of $\Delta(E)$.  By a final application of
Claim~\ref{theclaim}, we see that
$F^*_4(s,t)$ is a $\{2, 3\}$-unit.  Solving the implied Thue-Mahler equations,
one finds no unexpected solutions.

\subsection{The case $\bs n\bs=\bs2$}\label{neqtwo}

It remains to check that, for any $m$ and $Q$, the term $W_2$ has
a primitive divisor.  From~\eqref{vnform} we have, in the notation
of previous sections, that
$$\frac{U_2}{W_2}=\frac{36m\psi_2^3(Q)+\omega_2(Q)}{6\phi_2(Q)\psi_2(Q)},$$
where
$$\phi_2=x(x^3-8D),\quad \psi_2^2=4(x^3+D).$$
We suppose that every prime dividing $W_2$ also divides $W_1$, and
hence $6A_1B_1$.  If $p\geq 5$ is any prime with
$$\ord_p(\psi_2(Q))=\ord_p(4A_1^3+4DB_1^6)>0,$$ then $p$ is a
primitive divisor of $W_2$, or $p\mid 6m$.  Assuming that the
former is not the case, then, we see that $(4s+t)$ is a $\{2,
3\}$-unit, employing Claim~\ref{theclaim} as in the cases above.
That is, we know that $p\mid (4A_1^3+4DB_1^6)$ only if
$p\mid\gcd(A_1, m)$, in which case $\ord_p(s)>\ord_p(t)$.  In
fact, we can say slightly more.  Just as in
Section~\ref{thuesolve}, we may conclude from the inequality
$$(4A_1^3+4DB_1^6)=4C_1^2>0$$ that $(4s+t)$ is a positive $\{2,
3\}$-unit.

Now consider the primes dividing $A_1^3-8DB_1^6$.  Just as in the
previous case, if $p\geq 5$ is a prime divisor of this expression,
and $p$ is not a primitive divisor of $W_2$, then $p\mid m$.  If
this is the case, then another application of Claim~\ref{theclaim}
tells us that $(s-2t)$ is a $\{2, 3\}$-unit.  Solving  two linear
equations yields
\begin{alignat}{1}
9s&=W_1+2W_2\nonumber\\
9t&=-4W_1+W_2,\label{tea}
\end{alignat}
where $W_1$ and $W_2\in\Z$ are $\{2, 3\}$-units, and $W_2>0$.
It is~\eqref{tea} in which we are most
interested.  The second part of Claim~\ref{theclaim} tells us
that $\ord_p(t)$ is even for any $p\geq 5$.
Specifically, then, $t=-dx^2$, for some positive $d\mid 6$ and
some $x\in\Z$ (recall that $t<0$ as per
Section~\ref{thuesolve}).  Equation~\eqref{tea} is now the representation
of $-dx^2$ as a sum or
difference of $\{2, 3\}$-units.  We also know that $\gcd(W_1, W_2)\mid 9$,
as $\gcd(s, t)=1$.  We may
use this information to solve the above system of equations for all possible
values of $s$ and $t$.
First, a lemma.
\begin{lemma}
The only integral equations $a=b+c$ such that
\begin{enumerate}
\item $c<0$;
\item $b$ and $c$ are $\{2, 3\}$-units;
\item $\gcd(b, c)=1$; and
\item $a$, $2a$, $3a$, or $6a$ is a perfect square
\end{enumerate}
are the following:
\begin{alignat*}{2}
1^2 & = 2^1-1, \qquad & 7^2 &= 3^4-2^5,\\
1^2 & = 2^2-3^1, & 2.1^2 & = 3^1-1,\\
1^2 & = 3^1-2^1, & 2.2^2 & = 3^2-1,\\
1^2 & = 3^2-2^3, & 2.11^2 & = 3^5-1, \\
5^2 & = 3^3-2^1, & 3.1^2 & =2^2-1.\\
\end{alignat*}
\end{lemma}

\begin{proof}
Supposing that $a=b+c$ is one such equation, we may multiply both
sides by sufficient powers of 2 and 3
to obtain an equation of the form
$$q^2=r^3-2^\mu3^\nu,$$
where $r$ is a $\{2, 3\}$-unit.  Thus, if we write
$\mu=\mu_0+6\mu_1$, where $0\nobreak \leq \nobreak \mu_0 \nobreak
< \nobreak 6$, and write $\nu=\nu_0+6\nu_1$ similarly, we see that
$$\left(\frac{r}{2^{2\mu_1}3^{2\nu_1}}, \frac{q}{2^{3\mu_1}3^{3\nu_1}}\right)$$
is a $\{2, 3\}$-integral point on the elliptic curve
$Y^2=X^3-2^{\mu_0}3^{\nu_0}$, with the additional property that
$X$ is a $\{2, 3\}$-unit.  Using MAGMA~\cite{magma} to find all
$\{2, 3\}$-integral points on each curve of this form  (for $0\leq
\mu_0, \nu_0< 6$), and tracing these points back to the original
equations, we have our result.
\end{proof}

It is now a simple matter to solve the above system of equations for
all possible $s$ and $t$.
If $\gcd(W_1, W_2)=9$, then
$$-t=\frac{4}{9}W_1-\frac{1}{9}W_2$$
is an equation as in the lemma.  Moreover, as $\frac{4}{9}W_1$ is
divisible by 4 and $W_2>0$, most
of the equations listed in the lemma may be disqualified immediately.
We are left with the
possibilities $(W_1, W_2)=(9, 27)$ or $(9, 9)$.  These yield
$(s, t)=(7, -1)$ or $(3, -3)$, the latter
of which may be discarded as $\gcd(s,t)=1$ by construction.

If $\gcd(W_1, W_2)=3$, then
$$-3t=\frac{4}{3}W_1-\frac{1}{3}W_2$$
defines an equation as in the lemma, yielding again the solution
$(3, -3)$ for $(s,t)$.  The case $\gcd(W_1, W_2)=1$ yields no new
solution.

The only solution, then, that might contradict our claim, is that
corresponding to $(s, t)=(7, -1)$.  One may trace this solution
back to the point $(2, -1)$ on the curve $u^3+v^3=7$.  While the
terms $A_1^3+DB_1^6$ and $A_1^3-8DB_1^6$ do both turn out to be
$\{2, 3\}$-units in this example, we may simply compute $W_1=1$
and $W_2=3$ to see that $W_2$ does, indeed, have a primitive
divisor.

This completes the proof of Proposition~\ref{smallnprop}, and hence
the proof of Theorem~\ref{theanswer}.
\hfill\qed


\appendix
\section{Computations}\label{comps}

In this section we discuss some of the particulars of the
computations.

\subsection{Small values of $m$}\label{smallm}

As required in Sections~\ref{hereitcomes} and~\ref{hereitcomesZW}, we computed manually the
Zsigmondy bounds~$Z(A)$ and~$Z(W)$ for all the cases needed with $m\le 290$.
From the analysis given in Section~\ref{hereitcomes}, it is
sufficient to consider all $m\le 40$ (in fact we go out to $50$)
but only those $40\nobreak <\nobreak m<\nobreak 290$, for which
$m\equiv \pm 2\MOD 9$ and which possess a prime divisor
congruent to $1\MOD 6$. In all cases when the rank is positive the rational
torsion group is trivial.

For rank-1 curves, the values of $Z(W)$ are given in Figure~\ref{Zigfig1}.
The values of $Z(A)$ are much simpler: $Z(A)=0$ except in the case $m=7$, when $Z(A)=2$.

In rank-2 the situation is slightly more complicated. However,
provided $\hat h(Q)>0.1$, equation~\eqref{generic} yields the
Zsigmondy bound $Z(A)\nobreak \le \nobreak 12$, and
equation~\eqref{genericforV} yields the Zsigmondy bound $Z(W)\le
14$. In all cases, there were no non-trivial rational points with
$\hat h(Q)\le 0.1$ so no further checking was necessary. In
Figure~\ref{Zigfig2} we list generators for each of the rank-2
curves in the range. These were looked up in Cremona's
tables~\cite{cremona}, when the conductor was below $10^4$, or
computed using MAGMA~\cite{magma} for larger conductors.

\subsection{Division Polynomials}\label{divpolytables}

Various specific binary forms are used for computations in Section~\ref{formssection}, and these
are reproduced in Figures~\ref{fn} and~\ref{ftilde}.  To save space, only the coefficients have been
recorded.  The line
$$
\begin{array}{|c|c|}\hline
F & [v_d,  v_{d-1},\ \cdots\ , v_1, v_0]\\\hline \end{array}$$
in the table is to be interpreted as the statement
$$F(X,Y)=\sum_{i=0}^d v_iX^iY^{d-i}.$$
Note that, as a result of the action of complex multiplication on $E$, the polynomial $\psi_p$
factors when $p\equiv 1\MOD{3}$.  In actual computations, this fact can be exploited to great
advantage, allowing the corresponding Thue equations to be solved by elementary means.


\newpage
\begin{figure}\caption{Zsigmondy bounds in rank-1\label{Zigfig1}}
\begin{tabular}{|c|c|c||c|c|c|}\hline $m$& $E$-generator & $Z(W)$ & $m$& $E$-generator &$Z(W)$\\
\hline
6 &  [28,80] & 0 & 7 & [84,756] & 1\\
9 &  [36,108] & 1 & 12 &  [52,280] & 0 \\
13 &  [52,260] & 0 & 15 & [49,143] & 0\\
17 &  [84,684] & 0 &  20 & [84,648] & 0\\
22 &  [553/9,4085/27] & 0 & 26 & [156,1872] & 1\\
28  &  [84,504] & 1 & 31 &  [217,3131] & 0\\
33 &  [97,665] & 0 & 34 &  [2733,4455] & 0\\
35 &  [84,252] & 1 & 42 & [172,280] & 0\\
43 &  [129,1161] & 0 & 49 & [196,2548] & 0\\
50 &  [8148/27,138736/27] & 0 & 51 &   [5473/36,333935/216] & 0\\
58 & [9444/27,173600/27] & 0 & 61 & [732,19764] & 1 \\
79  &  [316,5372] & 0 & 97  & [388,7372] &  0 \\
133 & [228,2052] & 1 & 151 & [4228/9,261532/27] & 0\\
169 & [2028,91260]& 1 & 223 & [1561/4,49283/8] &   0 \\
241 & [6748,554300] & 1 & 259 & [777,20979] & 0 \\
277 & [5817/4,441261/8] & 0 & 286 & [588,12960] & 1\\
\hline
\end{tabular}
\end{figure}

\begin{figure}\caption{Rank-2 generators\label{Zigfig2}}
\begin{tabular}{|c|c|}\hline $m$& $E$-generators \\
\hline
19 & [156,1908], [228,3420] \\
37 & [84,36], [148,1628] \\
65 & [129,567], [156,1404] \\
124 & [372,6696],[2356,114328]\\
182 & [273,2457],[364,5824] \\
209 & [1596,63612],[532,11476] \\
218 & [1308,47088],[13881/25,1534221/125]\\
254 & [16257/4,2072385/8],[508,10160]\\
\hline
\end{tabular}
\end{figure}

\begin{figure}
\begin{tabular}{|l|l|}\hline
$F_5$ & [5, 95, -15, -25, -1]\\
$F_7$ &  [7,  986,  -2681, -12964, 3626, -1519, -686, -49, 1]\\
$F_8$ &   [2, 616, -7336, -1544, -3430, -4124, -952, -104, -1]\\
$F_{10}$ &  [1, 1173, -55284,29380,-368055,-1404072,-862941, 542232,\\
& \hfill\ldots -104805,-7070,-474,-177,1]\\
$F_{11}$ & [   11,  23221, -1153603, -62045313, 66133914, -1596123771,\\
& \hfill\ldots -8579472693, -4760052033, -22319781, 8054721004, \\
& \hfill\ldots10595519759, 4869514969, 1106263389, 189881835, \\
& \hfill\ldots59389374, 17393277, 2270301, 102729,605,-242, -1]\\
$F_{13}$ & [13, 74737, -10304874, -1459820466, 7383882519, -294761888811,\\
 &\hfill \ldots -3649379851026, -327751614216, 3634612800273, \\
  & \hfill\ldots 75587434125411, 206422282971957, 165623202699903,\\
   & \hfill\ldots 77423927253309, 50317031121903, 70684315657137, \\
    & \hfill\ldots 64207462488471, 30461492791431, 8167061938581, \\
     & \hfill\ldots 1237534488021, 33446767107, -47530886481, -16133119236, \\
      & \hfill\ldots -2480541102,-183218139, -6445998, -217503,\\
       & \hfill\ldots -22815,-338, 1]\\
$F_{14}$ &  [1,   8826,   -3182349,  27544616,  -1267563423,-29876807793,\\
 & \hfill\ldots -73452197357, -534368475927,   -321414204609,\\
  & \hfill\ldots    -159623734993,     -250499094747,  -930524257131,\\
   & \hfill\ldots -1172171589176,   -509647490898,  -20486729571,\\
    & \hfill\ldots  61406271479,   22270327506,  3403598121, 263510632,\\
     & \hfill\ldots   15278739, 2663808, 488510,  19851,  537, 1]\\
\hline
\end{tabular}\caption{The binary forms $F_n$\label{fn}}
\end{figure}

\begin{figure}
\begin{tabular}{|l|l|}\hline
$\tilde{F}_6$&[1, 57, 3, 1]\\
$\tilde{F}_9$&[1, 657, 6111, -3318, 19647, 12033, 3972, 684, 9, 1]\\
$\tilde{F}_{12}$ & [1, 3630, -28608, 392908, 212553, 1121508,168108, 62712, \\
 & \hfill\ldots 69507, 32782, 3684, 12, 1]\\
\hline
\end{tabular}\caption{The binary forms $\tilde{F}_n$\label{ftilde}}
\end{figure}
\end{document}